\documentclass[12pt]{amsart} 

\usepackage[margin=2.9cm]{geometry}             
\usepackage{color,ulem}
\usepackage{graphicx}
\usepackage{amssymb, enumerate,bm}
\usepackage[matrix,arrow,ps]{xy}
\usepackage{epstopdf, amsmath} 
\usepackage{paralist}
\newcommand{\C}{\mathbb C}

\newcommand{\R}{\mathbb R}
\newcommand{\B}{\mathbb B}
\newcommand{\transp}{\,^t}

\newcommand{\aaa}{\mathfrak{(a)}}

\newcommand{\bbb}{\mathfrak{(b)}}

\newcommand{\vnorm}[1]{{\| #1 \|}}

\DeclareMathOperator{\spanc}{span}
\newtheorem{theo}{Theorem}[section]
\newtheorem{lemma}[theo]{Lemma}
\newtheorem{cor}[theo]{Corollary}
\newtheorem{prop}[theo]{Proposition}

\theoremstyle{remark}
\newtheorem{remark}[theo]{Remark}
\theoremstyle{example}
\newtheorem{example}[theo]{Example}
\theoremstyle{definition}
\newtheorem{defi}[theo]{Definition}
\numberwithin{equation}{section}

\begin{document}

\begin{abstract}
We introduce a natural notion of defect for lifts of stationary discs attached to a CR nondegenerate submanifold  $M \subset \C^N$ of codimension $d$. 
We prove that the existence of a nondefective  lift of a stationary disc implies that germs of CR automorphisms of $M$ are uniquely determined by their $2$-jet. We emphasize that, in particular, this unifies all previously known results on this question.  
\end{abstract} 
\thanks{Research of the first author was supported by the Center for Advanced Mathematical Sciences.}

\author[Bertrand, Meylan]{Florian Bertrand and Francine Meylan}
\title[Nondefective  lifts of stationary discs and the 2-jet determination problem]{Nondefective  lifts of stationary discs and the 2-jet determination problem}

\subjclass[2010]{}

\thanks{}

\maketitle


\section*{Introduction}

The finite jet determination of CR automorphisms of a CR submanifold is an important problem in CR Geometry due to its connection with the equivalence problem. 
This question originates from the works of Cartan \cite{ca}, Chern-Moser \cite{ch-mo}, and Tanaka \cite{ta} who proved the $2$-jet determination for Levi nondegenerate hypersurfaces. Later on, the problem of finite jet determination in several different settings has attracted a lot of attention \cite{BER2, eb, eb-la-za, eb-la, KMZ1,za1,BER1, BMR, la-mi, ju, ju-la, ki-za, mi-za, la-mi2}; see also the surveys \cite{za2,la-mi3}. 
  
Recently, the method of stationary discs has been efficiently used to study the question of finite jet determination for finitely smooth submanifolds \cite{be-bl,be-de1,be-de-la,be-bl-me, tu3,be-me}. These invariant discs were introduced by Lempert \cite{le} (see also Tumanov \cite{tu}) as the complex geodesics of the Kobayashi metric in strongly convex domains. Using this approach, several papers address the $2$-jet determination for higher codimension submanifolds, leading to the best (to our knowledge) results in that direction. For instance, the $2$-jet determination has been obtained by Tumanov \cite{tu3} for strongly pseudoconvex submanifolds and by the authors \cite{be-me} in case the CR submanifold satisfies a certain nondegeneracy condition. In these  two papers, which "complement" each other,  the notion of {\it nondefective analytic discs} is essential. More precisely, 
in \cite{tu3}, Tumanov proves the  existence of a {\it nondefective stationary disc} attached to a strongly pseudoconvex submanifold which directly ensures the  $2$-jet determination of CR automorphisms. In general, the existence of 
such a disc does not guarantee the $2$-jet determination. For instance, Gregorovi\v{c} and Meylan \cite{gr-me} constructed a quadric of codimension  $4$ in $\C^{10}$ which does not satisfy 
the $2$-jet determination, yet admits a nondefective stationary disc (see Example 2.6 in \cite{be-me4}). Meanwhile, quadrics considered in \cite{be-me} admit a nondefective stationary disc and, along with their small perturbations, satisfy the $2$-jet determination.  

One of the aims of the present paper is precisely to understand the connection between such papers \cite{tu3, be-me, me, gr-me, be-me-ko}, and unify them through a common geometric condition.
The notion  of nondefective analytic discs goes back to the work of Tumanov \cite{tu0} where the existence of a {\it small} nondefective disc, which is ensured by the {\it minimality} of the CR sumbanifold, leads to  the wedge extendability of CR maps (see also \cite{ba-ro-tr} and Theorem 8.1.1 and Theorem 8.6.1 in \cite{ber}).  The main idea in the present paper is that since the existence of a nondefective stationary disc  does not guarantee the $2$-jet determination, we then seek further geometric information by focusing on the lifts of these discs. We then introduce a natural notion of defect for lifts  of stationary discs to the cotangent bundle. We observe that even if a stationary disc is nondefective, its lifts may or may not be nondefective. This point is crucial since it unifies results where $2$-jet determination holds \cite{tu3,be-me,be-me-ko} but also allows us to understand the reason why it fails in \cite{me, gr-me}:
\begin{theo}\label{theomain}
Let $M\subset \C^{N}$ be a $\mathcal{C}^4$ generic strongly Levi nondegenerate (at $0$) real submanifold. Assume that its model quadric $Q$ admits a nondefective stationary lift. Then germs at $0$ of CR automorphisms of $M$ of class $\mathcal{C}^3$ are uniquely determined by their $2$-jet at $0$.     
\end{theo}
This theorem implies both \cite{tu3} and \cite{be-me} (and also \cite{be-bl} in the hypersurface case). Indeed, we prove that if the CR submanifold $M$ is strongly pseudoconvex then a stationary disc is nondefective if and only if any of its lift is nondefective (see Corollary \ref{corstrong}). This crucial observation unveils the key point in \cite{tu3} that the nondefective disc which is constructed 
indeed admits a nondefective lift. We also prove that in case $M$ is nondegenerate in the sense of  \cite{be-me} then it admits a nondefective lift of a stationary disc (see Theorem \ref{theodefli}). Moreover, our result implies that the quadrics considered in \cite{me,gr-me} do not admit any nondefective lifts of a stationary disc. In addition, the assumption in Theorem \ref{theomain} is purely geometric and nontechnical, and does not imply any restriction on the codimension as in \cite{be-me} nor positivity condition as in \cite{tu3}. 

\vspace{0.5cm}

The paper is organized as follows. In the first section, we review several notions of nondegeneracy of a CR submanifold and discuss the theory of stationary discs. We emphasize in Remark \ref{rembig}, that discs we work with are not {\it small} unlike analytic discs considered in \cite{tu0,ba-ro-tr}. In the second section, we first recall the defect of a analytic disc and define  the corresponding notion for lifts of stationary discs. We then characterize lifts which are nondefective in Theorem 
\ref{theodefli} and Theorem \ref{theonondef}. The former result is the key in obtaining our main theorem which is stated in 2.5.

\section{Preliminaries}

\subsection{Nondegenerate generic submanifolds}
Let $M \subset \C^{n+d}$ be a real submanifold  of  codimension $d\ge 1$ given by the following $d$ real-valued defining functions
\begin{equation}\label{eqred}
\begin{cases}
r_1(z,w,\overline{z},\overline{w})=\Re e  w_1-\transp\bar z A_1 z+ \ldots\\
\ \ \ \ \vdots \\
r_d(z,w,\overline{z},\overline{w})=\Re e  w_d - \transp\bar z A_d z+ \ldots
\end{cases}
\end{equation}
The   $n\times n$ matrices  $A_1,\hdots,A_d$ are Hermitian and the dots represent terms in $z$ and $\Im m w$ of higher weighted order (see \cite{bo,bl-me}). We write $r:=(r_1,\ldots,r_d)$ and so $M=\{r=0\}$.  

In this section, we review  several key notions of nondegeneracy which are strongly connected with the unique jet determination of germs of CR automorphisms of $M$.

\begin{defi} [\cite{be1}]\label{defnondegbe}
A submanifold $M \subset \C^{n+d}$ given by \eqref{eqred} is  {\it Levi nondegenerate} at 0 if the following two conditions are satisfied
\begin{center}
\begin{tabular}{cl}
$\aaa$ & $A_1$,...,$A_d$ are linearly independent (equivalently on $\R$ or $\C$),\\
\\
$\bbb$ & $\bigcap_{j=1}^d\mathrm{Ker}A_j=\{0\}$.
\end{tabular}
\end{center}
\end{defi}

Note that Condition $\aaa$ implies  $d \le n^2$  and  that $M$ is  of finite type, while Condition $\bbb$ implies that $M$ is holomorphically nondegenerate. When $M$ is a hypersurface, Condition $\bbb$ implies Condition  $\aaa$. 

\begin{defi} [\cite{tu}] \label{defnondegtu}
A submanifold $M  \subset \C^{n+d}$ given by \eqref{eqred}  is  {\it strongly Levi nondegenerate} (resp. {\it strongly pseudoconvex}) at 0 if there exists a real linear combination of  
$A_1,\ldots,A_d$ which is invertible (resp. positive definite).
\end{defi}
Note that  when $M$ is strongly Levi nondegenerate then it satisfies Condition   $\bbb.$ We also point out that if the codimension of $M$  is maximal $d=n^2$ and if $M$ also satsifies Condition $\aaa$ then  it is strongly pseudoconvex (see in particular Example \ref{exd=n^2}).  
\begin{defi} [\cite{be-me}]
A submanifold $M  \subset \C^{n+d}$ given by \eqref{eqred}  is {\it$\mathfrak{D}$-nondegenerate} at 0 if it is strongly Levi 
nondegenerate at $0$ and  if there exists $V \in \C^n$  such that the $d\times d$ matrix 
$$\Re e (\transp \overline D A^{-1}D)$$ 
is invertible. Here, $A$ is an invertible real linear combination of $A_1,\ldots,A_d$, and 
$D$ is the $n \times d$ matrix whose columns $A_1V,A_2V,\ldots,A_dV$.
\end{defi}
We note that if $M$ is $\mathfrak{D}$-nondegenerate at $0$ then there exists $V \in \C^n$  such that the vectors $A_1V,\ldots,A_dV$ are $\R$-linearly independent in $\C^n$, which implies in particular the codimension restriction $d\leq 2n$. 

The last notion of nondegeneracy requires some setup  before being introduced. 
Consider  a  strongly Levi nondegenerate real submanifold $M\subset \C^{n+d}$ given by \eqref{eqred}. We fix $b \in \R^d$ such that $\sum_{j=1}^d b_jA_j$  is invertible. For $a \in \C^d$ sufficiently small, we set 
$$P:=\sum_{j=1}^d a_j A_j  \ \mbox { and }  \ A:=\sum_{j=1}^d (b_j -a_j-\overline{a_j})A_j.$$ Since the norm $\|P\|$ is sufficiently small, the  following  quadratic matrix equation 
\begin{equation}\label{eqX}
PX^2 + AX +   \transp{\overline {P}}=0
\end{equation} 
has a unique $n\times n$ matrix solution $X$  satisfying  $\|X\|<1.$  Note that the solution $X$ depends on both $a$ and $b$. 
  The following definition, inspired by the work of Tumanov in the strongly pseudoconvex framework  \cite{tu}, was introduced in \cite{be-me}.
 \begin{defi}[\cite{be-me}]\label{defstatmin}
 Let  $M\subset \C^{n+d}$ be a  strongly Levi nondegenerate real submanifold given by \eqref{eqred} and let $V\in \C^n$. 
 In the above context, $M$ is {\it stationary minimal at $0$ for $(a,b-a-\overline{a},V)$} if the matrices 
 $A_1,\ldots,A_d$ restricted to the  orbit space
 $${\rm span}_\R \{V, XV, {X^2}V, \ldots, {X^k}V, \ldots\}\subset \C^n$$ are 
 $\R$-linearly independent. 
 \end{defi}
 Note that $M$ is stationary minimal for $(0,b,V)$ if and only if it is $\mathfrak{D}$-nondegenerate.
Tumanov proved in Theorem 5.3 \cite{tu3} that if  $M$ satisfies Condition $\aaa$, then there exist $a\in \Bbb R^d $ and $V \in \C^n$ such that  the matrices 
 $A_1,\ldots,A_d$ restricted to the orbit space
$${\rm span}_\R \{V, PV, {P^2}V, \ldots, {P^k}V, \ldots\}$$ are 
 $\R$-linearly independent. In particular, this important result implies that strongly pseudoconvex submanifolds  which satisfy Condition $\aaa$ are stationary minimal. 

\vspace{0.5cm}
To conclude, we emphasize that all above definitions are independent of the choice of holomorphic coordinates.

 \subsection{Stationary discs}
We denote by $\Delta$ the unit disc in $\C$ and by $\partial \Delta$ its boundary. 
 For an integer $k \geq 0$ and a real number $0< \alpha<1$, we denote by $\mathcal C^{k,\alpha}$ the space of real-valued functions  defined on $\partial\Delta$ of class $\mathcal{C}^{k,\alpha}$. This space is equipped with  its usual norm
$$\|f\|_{\mathcal{C}^{k,\alpha}}=\sum_{j=0}^{k}\|f^{(j)}\|_\infty+
\underset{\zeta\not=\eta\in \partial\Delta}{\mathrm{sup}}\frac{\|f^{(k)}(\zeta)-f^{(k)}(\eta)\|}{|\zeta-\eta|^\alpha},$$
where $\|f^{(j)}\|_\infty=\underset{\zeta \in \partial\Delta}{\mathrm{max}}|f^{(j)}(\zeta)|$.
We write $\mathcal C_\C^{k,\alpha} = \mathcal C^{k,\alpha} + i\mathcal C^{k,\alpha}$ and equip this space with   the following norm  
\begin{equation*}
\|f\|_{\mathcal{C}_{\C}^{k,\alpha}}=
\|\Re e  f\|_{\mathcal{C}^{k,\alpha}}+\|\Im m f\|_{\mathcal{C}^{k,\alpha}}.
\end{equation*}
We denote by $\mathcal A^{k,\alpha} \subset \mathcal C_{\C}^{k,\alpha}$ the subspace of {\it analytic discs} that consists of functions $\mathcal C_\C^{k,\alpha}$ which extend holomorphically to $\Delta$. Alternatively, analytic discs are functions defined on $\overline{\Delta}$ which are holomorphic on $\Delta$ and whose restriction to $\partial\Delta$ belongs to $\mathcal C_\C^{k,\alpha}$. In any case, we write $f:\Delta \to \C$, omitting the closure notation $\overline{\Delta}$,  whenever we deal with an analytic disc.

\vspace{0.5cm}

Consider now a real submanifold $M=\{r=0\}\subset \C^{n+d}$ given by  (\ref{eqred}). An analytic disc $f: \Delta \to \C^{n+d}$ 
is {\it attached to  $M$} whenever $f(\partial \Delta) \subset M$. Following Lempert \cite{le} and Tumanov \cite{tu}, we recall the definition of stationary discs. 
\begin{defi}\label{defistat}
An analytic disc $f: \Delta \to \C^{n+d}$ attached to  $M$ is a {\it stationary disc} if  it admits an 
analytic lift $\bm{f}=(f,\tilde{f}): \Delta \to T^*\C^{n+d}$ to the cotangent bundle, continuous up to 
 $\partial \Delta$, and such that for all $\zeta \in \partial\Delta$ we have $\bm{f}(\zeta) \in \mathcal{N}M(\zeta)$, where 
\begin{equation*}
\mathcal{N}M(\zeta):= \left\{(z,w,\tilde{z},\tilde{w}) \in T^*\C^{n+d} \ | \ (z,w) \in M, (\tilde{z},\tilde{w}) \in 
\zeta N^*_{(z,w)} M\setminus \{0\}\right\},
\end{equation*}
 and where
$$N^*_{(z,w)} M:=\spanc_{\R}\{\partial r_1(z,w,\overline{z},\overline{w}), \ldots, \partial r_d(z,w,\overline{z},\overline{w})\}$$ is the conormal fiber at $(z,w)$ of $M$. 
We call the map $\bm{f}=(f,\tilde{f})$ a {\it stationary lift for $M$} and we denote by $\mathcal{S}(M)$ the set of stationary lifts with $f$ nonconstant. 
\end{defi}
It is important to note that an analytic disc $f$ attached to  $M$ is stationary  if and only if there exist $d$ real-valued functions $c_1, \ldots, c_d : \partial \Delta \to \R$ satisfying $\sum_{j=1}^dc_j(\zeta)\partial r_j(0)\neq 0$ for all $\zeta \in \partial \Delta$   and such that the map 
$$\zeta \mapsto \zeta\sum_{j=1}^dc_j(\zeta)\partial r_j\left(f(\zeta), \overline{f(\zeta)}\right)$$
defined on $\partial \Delta$ extends holomorphically to $\Delta$.

\vspace{0.5cm}

Consider now a {\it quadric} $Q \subset \C^{n+d}$ of codimension $d$, that is, a submanifold given by  
\begin{equation}\label{eqquad}
\begin{cases}
\rho_1(z,w,\overline{z},\overline{w})=\Re e  w_1-\transp\bar z A_1 z\\
\ \ \ \ \vdots \\
\rho_d(z,w,\overline{z},\overline{w})=\Re e  w_d - \transp\bar z A_d z
\end{cases}
\end{equation} 
where $A_1,\ldots,A_d$ are $n\times n$ Hermitian matrices. We set $\rho:=(\rho_1,\ldots,\rho_d)$ and write $Q=\{\rho=0\}$. Assume that $Q$ is strongly Levi nondegenerate 
and fix $b \in  \Bbb R^d$  such that $\sum_{j=1}^d b_jA_j$ is invertible. We consider the following subset of stationary lifts
$$\mathcal{S}_0(Q) = \{\bm{f} \in  \mathcal{S}(Q) \ | \ \bm{f}(1)=(0,0,0,b/2)\}.$$  
In this context, it turns out that the following analytic disc  
\begin{equation}\label{eqinitlift}\small
 \bm {f_0}(\zeta)=\left((1-\zeta)V_0,2(1-\zeta) \transp\overline{V_0}A_1V_0,\ldots,2(1-\zeta) \transp \overline{V}A_dV_0,(1-\zeta)  \transp \overline{V_0}\sum_{j=1}^d b_jA_j, \frac{\zeta}{2} b\right)
 \end{equation} 
with $V_0 \in\C^n$, is a stationary lift in $\mathcal{S}_0(Q)$.
 In \cite{be-me2}, we proved that stationary lifts $\bm{f}=(h,g,\tilde{h},\tilde{g}) \in \mathcal{S}_0(Q)$ which are close to  $\bm {f_0}$ are exactly   of the form
\begin{equation}\label{eqform}
\begin{cases}
h(\zeta)=(1-\zeta)(I-\zeta X)^{-1}V=(1-\zeta)\sum_{r=0}^{\infty}\zeta^rX^r V\\
\\
g_j(\zeta)={ \transp\overline{V}}A_j V-2 \transp \overline {V}A_j \zeta (I-\zeta X)^{-1}(I-X)V + \transp \overline {V}( \transp\overline {X}K_j-K_jX))V +\\
\hspace{1.4cm} \transp \overline {V}(I- \transp \overline {X})K_j(I + 2\zeta X(I-\zeta X)^{-1})(I-X) V\\
\\
\tilde{h}(\zeta)=-\transp\overline{h(\zeta)} \left(P+A\zeta+\transp \overline{P}\zeta^2\right)\\
\\
\displaystyle \tilde{g}(\zeta)=\frac{a+(b-a-\overline a)\zeta+\overline{a}\zeta^2}{2}\\
\end{cases}
\end{equation}
with  $V \in \Bbb C^n$ close to $V_0$ and $a \in \C^d$ is sufficiently small, where $X$ is the unique $n\times n$ matrix solution of \eqref{eqX} (with $b- a-\overline{a}$) satisfying  $\|X\|<1$, and where 
\begin{equation}\label{eqK_j}
K_j=\sum_{r=0}^\infty {\transp \overline {X}}^r A_j X^r, \ j=1, \dots,d.  
\end{equation}

\vspace{0.5cm}

Stationarity may be rephrased into a nonlinear boundary value problem. This is essentially due to the fact that the fibration $\mathcal{N}M(\cdot)$  introduced in Definition \ref{defistat} can be expressed by means of $2n+2d$ real defining equations. More precisely, in the case of  a quadric $Q= \{\rho=0\}\subset \C^{n+d}$ given by \eqref{eqquad}, the defining functions for the fibration $\mathcal{N}Q(\cdot)$ are given by:  
 \begin{equation*}\small
\left\{
\begin{array}{lll} 
\displaystyle \tilde{\rho}_1(\zeta)(\bold{Z},\overline{\bold {Z}}) & = & \rho_1(z,w,\overline{z},\overline{w})\\
\vdots& \vdots&\vdots\\
\tilde{\rho}_d(\zeta)(\bold{Z},\overline{\bold {Z}})  & = & \rho_d(z,w,\overline{z},\overline{w})\\
\tilde{\rho}_{d+1}(\zeta)(\bold{Z},\overline{\bold {Z}})  & = & \left(\tilde{z_1}+2 \transp\bar{z} \sum_{j=1}^d \tilde{w_j}(A_j)_1\right)+\left(\overline{\tilde{z_1}}+2 \transp z \sum_{j=1}^d \overline{\tilde{w_j}(A_j)_1}\right)\\
\vdots& \vdots&\vdots\\
\tilde{\rho}_{d+n}(\zeta)(\bold{Z},\overline{\bold {Z}})  & = & \left(\tilde{z_n}+2 \transp\bar{z} \sum_{j=1}^d \tilde{w_j}(A_j)_n\right)+\left(\overline{\tilde{z_n}}+2 \transp z \sum_{j=1}^d \overline{\tilde{w_j}(A_j)_n}\right)\\
\tilde{\rho}_{d+n+1}(\zeta)(\bold{Z},\overline{\bold {Z}})  & = &  i\left(\tilde{z_1}+2 \transp\bar{z} \sum_{j=1}^d \tilde{w_j}(A_j)_1\right)-i\left(\overline{\tilde{z_1}}+2 \transp z \sum_{j=1}^d \overline{\tilde{w_j}(A_j)_1}\right)\\
\vdots& \vdots&\vdots\\
\tilde{\rho}_{2n+d}(\zeta)(\bold{Z},\overline{\bold {Z}})  & = &  i\left(\tilde{z_n}+2 \transp\bar{z} \sum_{j=1}^d \tilde{w_j}(A_j)_n\right)-i\left(\overline{\tilde{z_n}}+2 \transp z \sum_{j=1}^d \overline{\tilde{w_j}(A_j)_n}\right)\\
\tilde{\rho}_{2n+d+1}(\zeta)(\bold{Z},\overline{\bold {Z}})  & = & \displaystyle i\overline{\zeta}\tilde{w_1}-i\zeta\overline{\tilde{w_1}}\\
\vdots& \vdots&\vdots\\
\tilde{\rho}_{2n+2d}(\zeta)(\bold{Z},\overline{\bold {Z}})  & = & \displaystyle i\overline{\zeta}\tilde{w_d}-i\zeta\overline{\tilde{w_d}}\\
\end{array}
\right.
\end{equation*}
Here, $\bold{Z}=(z,w,\tilde{z},\tilde{w})$ and $(A_j)_l$ denotes, for $j=1,\dots,d$ and $l=1,\ldots,n$, the $l^{\rm th}$ column  of $A_j$.
We set $\tilde{\rho}:=(\tilde{\rho}_1,\cdots,\tilde{\rho}_{2n+2d})$ and we have $(z,w,\tilde{z},\tilde{w}) \in \mathcal{N}Q(\zeta)$ if and only if $\tilde{\rho}(\zeta)(z,w,\tilde{z},\tilde{w})=0$.  
  For a general submanifold $M=\{r=0\}\subset \C^{n+d}$ given by \eqref{eqred},  we denote by $\tilde{r}$ the corresponding defining functions of $\mathcal{N}M(\cdot)$. 
It follows that an analytic disc $\bm{f}:\Delta \mapsto T^*\C^{n+d}$ is a stationary 
lift for $M$ if and only if it satisfies
\begin{equation}\label{eqstatrh}
\tilde{r}(\zeta)\left(\bm{f}(\zeta),\overline{\bm{f}(\zeta)}\right)=0
\end{equation}
for all $\zeta \in \partial \Delta.$ The problem \eqref{eqstatrh} is a nonlinear Riemann-Hilbert problem and was settled in \cite{be-bl-me} for strongly Levi nondegenerate submanifolds, leading to the following construction of stationary lifts:
\begin{theo}[\cite{be-bl-me}]\label{theodiscscons}
Let $Q=\{\rho=0\}\subset \C^{n+d}$ be a strongly Levi nondegenerate quadric  of real codimension $d$ given by (\ref{eqquad}).  Consider an initial stationary lift 
$\bm{f_0}=(h_0,g_0,\tilde{h_0},\tilde{g_0}) \in \mathcal{S}_0(Q)$ of the form  (\ref{eqinitlift}). Then there exist  an open 
neighborhood $U$ of $\rho$ in $(\mathcal{C}^4(\B))^{d}$ and a real number $\varepsilon>0$ such that for any $r \in U$, the set 
$$ \{\bm{f} \in \mathcal{S}_0(\{r=0\})\ \ | \  
\|\bm{f}-\bm{f_0}\|_{\mathcal A^{1,\alpha}_{0}}<\varepsilon\}$$
is a $\mathcal{C}^1$ real submanifold of dimension $2n+2d$ of the Banach space of analytic discs.    
\end{theo}
In the above theorem, $\mathcal A^{1,\alpha}_{0}$ denotes the affine space of analytic discs  $\bm{f} \in  (\mathcal{A}^{1,\alpha})^{2n+2d}$  of the form
$$\bm{f}(\zeta)=\left((1-\zeta)h,(1-\zeta)g,(1-\zeta)\tilde{h},(1-\zeta)\tilde{g}+\frac{\zeta}{2} c\right)$$
equipped with the norm $\left\|\bm{f}\right\|_{\mathcal A^{1,\alpha}_{0}}=\vnorm{(h,g,\tilde{h},\tilde{g})}_{\mathcal{C}_{\C}^{1,\alpha}}.$
We end this section with a crucial remark. 
\begin{remark}\label{remman}
In case the initial disc $\bm{f_0}$ is of the more general form \eqref{eqform} with $a \in \C^d$ sufficiently small, then  the set 
$$ \{\bm{f} \in \mathcal{S}_0(\{r=0\})\ \ | \  \|\bm{f}-\bm{f_0}\|_{\mathcal A^{1,\alpha}_{0}}<\varepsilon\}$$
is also a $\mathcal{C}^1$ real submanifold of dimension $2n+2d$ provided that the defining function $r$ is close enough to $\rho$. 
\end{remark}

\begin{remark}\label{rembig}
We recall that  a submanifold $M=\{r=0\}$ given by \eqref{eqred} is strongly Levi nondegenerate if and only if its conormal bundle  $N^*M\subset T^*\C^{n+d}$ is totally real at
 $\left(0,\sum_{j=1}
^db_j\partial r_j(0)\right)$, where $b\in \R^d$ is such that $\sum_{j=1}^d b_jA_j$  is invertible (see Proposition 1.1 \cite{tu}). In such case, the  fibration $\mathcal{N}M(\cdot)$ 
is such that  $\mathcal{N}M(\zeta)\subset T^*\C^{n+d}$ is totally real for every $\zeta \in \partial \Delta$, which is equivalent to the fact that   $\cup_{\zeta \in \partial \Delta} 
\mathcal{N}M(\zeta) \times \partial \Delta \subset T^*\C^{n+d} \times \C$ is a totally real submanifold. 
We know that there are no small analytic discs attached to a totally real submanifold due to its local polynomial convexity. Discs constructed in Theorem 
\ref{theodiscscons} are relatively large and their construction depends on  techniques, developed by Forstneri\v{c} \cite{fo} and Globevnik \cite{gl1,gl2}, which rely on 
the values of certain integers, namely the partial indices and the Maslov index, associated with
 the initial disc ${\bm f_0}$.  
\end{remark}

\section{Nondefective stationary lifts}

In this section, we review the notions of defective and nondefective analytic discs and introduce the corresponding notions for 
stationary lifts. We then characterize nondefective stationary lifts and use this characterization to obtain new results on the  $2$-jet determination of CR automorphisms.  
\subsection{Nondefective analytic discs}
Consider a real submanifold $M=\{r=0\}\subset \C^{n+d}$ of the form \eqref{eqred}. 
According to \cite{tu0} (see also \cite{ber, ba-ro-tr}), an analytic disc $f:\Delta \to \C^{n+d}$ attached to $M$ is {\it defective} if there exists a continuous 
map 
$c: \partial \Delta \to \R^{d}$ such that the map 
$$\zeta \mapsto c(\zeta)\partial r\left(f(\zeta),\overline{f(\zeta)}\right)=\sum_{j=1}^d c_j(\zeta) \partial r_j\left(f(\zeta),\overline{f(\zeta)}\right)$$
defined on $\partial \Delta$ extends holomorphically to $\Delta$. The disc is  {\it nondefective} if it is not defective. We recall that the existence of a small nondefective disc ensures the wedge extendability of CR maps  \cite{tu0} (see also Theorem 8.6.1 \cite{ber}).   
The {\it defect} of an analytic  disc attached to $M$ is then defined by
$${\rm def}(f):={\rm dim}_\R \left\{c:\partial \Delta \to \R^d \ | \ \zeta \mapsto c(\zeta)\partial r\left(f(\zeta),\overline{f(\zeta)}\right) \mbox{ extends holomorphically to } \Delta\right\}.$$
If follows that  $f$ is nondefective if and only if  ${\rm def}(f)=0$. It is important to point out the defect of an analytic disc is a biholomorphically invariant integer.  
Moreover, note that a defective analytic disc is necessarily a stationary disc. 
In case $Q\subset \C^{n+d}$ is a  quadric, we have the following lemma:
\begin{lemma}
Let 	$Q=\{\rho=0\} \subset \C^{n+d}$  be a quadric given by \eqref{eqquad}.
An analytic disc $f=(h,g)$ attached to  $Q$  is defective if and only if there exists a constant $c \in \R^d\setminus\{0\}$ such that the map 
$$\zeta \mapsto c\partial_z \rho\left(h(\zeta),\overline{h(\zeta)}\right)= \sum_{j=1}^d c_j \partial_z \rho_j\left(h(\zeta),\overline{h(\zeta)}\right)=-\transp\overline{h(\zeta)}{\sum_{j=1}^d c_jA_j}$$
defined on $\partial \Delta$ extends holomorphically to $\Delta$.
\end{lemma}
\begin{proof}
This indeed follows from the fact that $\partial_w \rho_j, j=1,\ldots,d$,  is a constant vector since $Q$ is a quadric.  
\end{proof} 
As a consequence, the defect of an analytic disc  attached to  a quadric  $Q=\{\rho=0\}$ is
$${\rm def}(f)={\rm dim}_\R \left\{c \in \R^d \ | \ \zeta \mapsto \transp\overline{h(\zeta)}{\sum_{j=1}^d c_jA_j} \mbox{ extends holomorphically to } \Delta\right\}.$$

We then recall the following  lemma which characterizes the defect of analytic discs of a special given form. 
\begin{lemma}[\cite{be-me,be-me4}]\label{lemnondef}
Let $Q \subset \C^{n+d}$ be a quadric given by \eqref{eqquad} and let $f=(h,g)$ be an analytic disc attached to $Q$ with $h$ of the form 
\begin{equation}\label{eqdisc}
h(\zeta)=(1-\zeta)V
\end{equation}
 for some $V \in \C^n$. 
 \begin{enumerate}[i.]
 \item  The disc $f$ is nondefective if and only if $A_1V,\ldots,A_dV$ are $\R$-linearly independent.
 \item We have 
 $${\rm def}(f)=\dim {\rm ker}D_{|\R^d}=d-{\rm dim} \ {\rm span}_\R\{A_1V,\ldots,A_dV\},$$
where  $D$ is the $n\times d$ matrix whose $j^{th}$ column is $A_jV$, $j=1,\ldots,d$. 
 \item The rank of the real differential at $V$ of the Levi map $L_0: \C^n \to \R^d$ defined by 
 $$L_0(W)=\left(\transp \overline{W}A_1W,\transp \overline{W}A_2W,\ldots,\transp \overline{W}A_dW\right)$$  is  
$${\rm rank}(d_VL_0)=d - {\rm def}(f).$$ 
In particular, the range of the Levi map of $Q$ contains an open set if and only if $Q$ admits a nondefective analytic disc of the form \eqref{eqdisc}.
 \end{enumerate}
\end{lemma} 
Note that analytic discs  $f=(h,g)$ with $h$ of the form \eqref{eqdisc} are stationary discs for $Q$ with a lift given exactly by \eqref{eqinitlift}. For a stationary disc $f$ (with a lift) of the more general form \eqref{eqform}, the corresponding algebraic characterization was obtain in Lemma 2.6 and Proposition 3.3 \cite{be-me2}:  
\begin{lemma}[\cite{be-me2}]\label{lemstatmin}
Let $Q \subset \C^{n+d}$ be a strongly Levi nondegenerate quadric given by \eqref{eqquad} and let $f=(h,g)$ be a stationary disc for $Q$ of the form \eqref{eqform}. 
The disc $f$ is nondefective if and only if $Q$ is stationary minimal at $0$  for $(a,b-a-\overline{a},h(0))$. Moreover we have 
$${\rm def}(f)=\dim {\rm ker} \ \Re e \left (\sum_{r=0}^{\infty} \transp \overline V{\transp \overline {X}}^rA_jA_sX^r V\right)_{j,s}.$$
\end{lemma}

\subsection{Nondefective stationary lifts}
Consider  a real submanifold $M=\{r=0\}\subset \C^{n+d}$ of the form \eqref{eqred}. 
We now introduce a natural notion of defect for stationary lifts. Recall  that such a lift $\bm{f} \in \mathcal{S}(M)$ is "attached" to the fibration $\mathcal{N}M(\cdot)$ in the sense that 
$\bm{f}(\zeta) \in \mathcal{N}M(\zeta)$ for all $\zeta \in \partial \Delta$. Recall also that the $2n+2d$ defining functions of $\mathcal{N}M(\cdot)$ are denoted by $\tilde{r}.$  In this context, we introduce the following more general notion of defect.
\begin{defi}\label{defidef} 
Let $M=\{r=0\}\subset \C^{n+d}$ be a  real submanifold of the form \eqref{eqred}. 
A stationary lift $\bm{f} \in \mathcal{S}(M)$  is {\it defective} if there exists a continuous map $\tilde{c}: \partial \Delta \to \R^{2n+2d}$ such that the map  
$$\zeta \mapsto \tilde{c}(\zeta)\partial \tilde{r}\left(\bm{f},\overline{\bm{f}}\right)(\zeta)= \tilde{c}(\zeta)\partial \tilde{r}(\zeta)\left(\bm{f}(\zeta),\overline{\bm{f}(\zeta)}\right)$$
defined on $\partial \Delta$ extends holomorphically to $\Delta$. The {\it defect} of $\bm{f}$, denoted by $\bm{{\rm def}}(\bm{f})$, is defined by
$${\rm dim}_\R \{\tilde{c}:\partial \Delta \to \R^{2n+2d} \ | \ \zeta \mapsto \tilde{c}(\zeta)\partial \tilde{r}\left(\bm{f},\overline{\bm{f}}\right)(\zeta) \mbox{ extends holomorphically to } \Delta\}.$$
\end{defi}

We write  $\tilde{c}=(\gamma,u,v,\delta)$ with $\gamma,\delta: \partial \Delta \to \R^d$ and $u,v: \partial \Delta \to \R^n$. 
We  first observe the following basic connection between the defect of a stationary disc and of its lifts.  
\begin{lemma}
If $f$ is a defective stationary disc for a submanifold $M\subset \C^{n+d}$ of the form \eqref{eqred} then any of its lifts $\bm{f}=(f,\tilde{f})$ is defective. Moreover, we have 
$${\rm def}(f) \leq \bm{{\rm def}}({\bm f}).$$ 
\end{lemma}
\begin{proof}
Assume that $f$ is a defective stationary disc for 
$M=\{r=0\}\subset \C^{n+d}$  and let $c: \partial \Delta \to \R^{d}$ be a continuous map such that $\zeta \mapsto c(\zeta)\partial r(f(\zeta),\overline{f(\zeta)})$
extends holomorphically to $\Delta$. Note that the first $d$ equations of $\tilde{r}$ are exactly $r$ and do not involve any of the fiber variables $\tilde{z}$ and $\tilde{w}$. It follows that for $\tilde{c}=(c,0,0,0)$, the map 
$$\zeta \mapsto \tilde{c}(\zeta)\partial \tilde{r}\left(\bm{f},\overline{\bm{f}}\right)(\zeta)=\left(c(\zeta)\partial r\left(f(\zeta),\overline{f(\zeta)}\right),0\right)$$
extends holomorphically to $\Delta$. This also implies, due to dimensional reasons, that ${\rm def}(f) \leq \bm{{\rm def}}({\bm f})$.
\end{proof} 
We now  provide an example of a quadric which admits a nondefective stationary disc $f$ with a defective lift ${\bm f}$, that is, with 
$0={\rm def}(f)<\bm{{\rm def}}(\bm{f}).$
\begin{example}\label{exline}
Consider the  quadric of codimension $3$ in $\Bbb C^7$ defined by  
\begin{equation*}
\begin{cases}
\rho_1=\Re e  w_1-|z_1|^2\\
\rho_2=\Re e  w_2-|z_2|^2\\
\rho_3=\Re e  w_3 - z_1\overline{z_4} - \overline{z_1}z_4 - z_2\overline{z_3} - \overline{z_3}z_2
\end{cases}
\end{equation*} 
 For $V=\transp(1,1,0,0)$, the  stationary disc $f=((1-\zeta) V,g)$ is nondefective since $A_1V, A_2V,A_3V$ are $\R$-linearly independent. Consider now the stationary lift $\bm{f}=((1-\zeta) V,g,\tilde{h},\tilde{g}),$ 
  with $\tilde{g}(\zeta)=(0,0,\zeta/2)$. In that case, we have 
  $$\partial{\tilde {\rho_1}} \left(\bm{f}, \bar {\bm{f}}\right)(\zeta)= (\underbrace{-1+ \bar\zeta,0,0,0}_{\partial_z{\tilde {\rho_1}}},\underbrace{1/2, 0, 0}_{\partial_w{\tilde {\rho_1}}},\underbrace{0,0,0,0}_{\partial_{\tilde{z}}{\tilde {\rho_1}}},\underbrace{0,0,0}_{\partial_{\tilde{w}}{\tilde {\rho_1}}})$$
 $$\partial{  \tilde {\rho_2}} \left(\bm{f}, \bar {\bm{f}}\right)(\zeta)= (0,-1+ \bar\zeta,0,0,0,1/2, 0,0,0,0,0,0,0,0)$$
 $$\partial{  \tilde {\rho_6}} \left(\bm{f}, \bar {\bm{f}}\right)(\zeta)= (0, \overline{\zeta}, 0, 0,  0,0,0,0,0,1,0,0,0,2- 2\bar\zeta)$$
  $$\partial{  \tilde {\rho_7}} \left(\bm{f}, \bar {\bm{f}}\right)(\zeta)= (\overline{\zeta},0, 0, 0, 0,0,0,0,0,0,1,0,0,2- 2\bar\zeta).$$
  Therefore, for 
$$\tilde{c}=(1,-1,0,0,0,1, -1,0,0,0,0,0,0,0),$$ 
the map 
$$\zeta \mapsto \tilde{c}\partial \tilde{\rho}(\zeta)(\bm{f}(\zeta),\overline{\bm{f}(\zeta)})= \partial{  \tilde {\rho_1}} (\bm{f}, \bar {\bm{f}})(\zeta) - \partial{  \tilde {\rho_2}} (\bm{f}, \bar {\bm{f}})(\zeta)+\partial{  \tilde {\rho_6}} (\bm{f}, \bar {\bm{f}})(\zeta)-\partial{  \tilde {\rho_7}} (\bm{f}, \bar {\bm{f}})(\zeta)$$ 
extends holomorphically to $\Delta$. Thus, the stationary lift $\bm{f}$ is defective.  To be more precise, we will see in Example \ref{exdef} that $\bm{{\rm def}}(\bm{f})=1$. 
\end{example}

In the following proposition, we establish another basic connection with the standard notion of defect. 
\begin{prop}\label{proplift}
Let $M=\{r=0\}\subset \C^{n+d}$ be a submanifold given by \eqref{eqred}. If a stationary lift $\bm{f} \in \mathcal{S}(M)$ is nondefective then the analytic disc 
$\zeta \mapsto (\bm{f}(\zeta),\zeta) \in T^*\C^{n+d} \times \C$ attached to the $\cup_{\zeta \in \partial \Delta} \mathcal{N}M(\zeta) \times \partial \Delta$ is 
nondefective (in the standard sense).
\end{prop}
\begin{proof}
We first note that $\cup_{\zeta \in \partial \Delta} \mathcal{N}M(\zeta) \times \partial \Delta \subset T^*\C^{n+d} \times \C$ is described by the $2(n+d)$ real-valued functions $\tilde{r}$ and the additional real-valued function
$$\tilde{r}_{2n+2d+1}(Z,\zeta,\overline{Z},\overline{\zeta}):=\zeta \overline{\zeta}-1=0.$$ In other words, 
$(Z,\zeta) \in \cup_{\zeta \in \partial \Delta} \mathcal{N}M(\zeta) \times \partial \Delta$ if and only if it satisfies 
\begin{equation*}
\begin{cases}
\tilde{r}(Z,\zeta,\overline{Z},\overline{\zeta})=0 \\
\tilde{r}_{2n+2d+1}(Z,\zeta,\overline{Z},\overline{\zeta})=0
\end{cases}
\end{equation*}
Here, we write 
$\tilde{r}(Z,\zeta,\overline{Z},\overline{\zeta})$ instead of  $\tilde{r}(\zeta)(Z,\overline{Z})$ to emphasize that $\zeta$ is no longer a parameter. We then set $\tilde{\bm{r}}:=(\tilde{r},\tilde{r}_{2n+2d+1})$ and we note that $\partial\tilde{\bm{r}}$ is the following $(2n+2d+1)\times (2n+2d+1)$ matrix 
 $$\partial\tilde{\bm{r}}=
 \begin{pmatrix}
 \partial \tilde{r} &  \partial_\zeta \tilde{r}\\
 0 & \overline{\zeta}
 \end{pmatrix}$$

Assume now that $\bm{f} \in \mathcal{S}(M)$ is such that the analytic disc $\zeta \mapsto  (\bm{f}(\zeta),\zeta)$ is defective and let 
$$(\tilde{c},\lambda)=(\gamma,u,v,\delta,\lambda): \partial \Delta \to \R^{2n+2d+1}$$ be such that the map 
$$\zeta \mapsto (\tilde{c}(\zeta),\lambda(\zeta))\partial \tilde{\bm r}\left((\bm{f}(\zeta),\zeta),\left(\overline{\bm{f}}(\zeta),\overline{\zeta}\right)\right)=\left(\tilde{c}(\zeta)\partial \tilde{r}\left(\bm{f},\overline{\bm{f}}\right)(\zeta),\tilde{c}(\zeta) \partial_\zeta \tilde{r}\left(\bm{f},\overline{\bm{f}}\right)(\zeta)+\lambda(\zeta)\overline{\zeta}\right)$$
 extends holomorphically to $\Delta$. It follows that $\bm{f}$ is defective 
since, in particular, the map
$\zeta \mapsto \tilde{c}(\zeta)\partial \tilde{r}(\bm{f},\overline{\bm{f}})(\zeta)$ extends holomorphically to $\Delta$.

\end{proof}

After Theorem \ref{theodefli}, we discuss in Remark \ref{remconv}   the converse of this result in the case of a quadric.

\vspace{0.5cm}

In the next two sections, we aim to characterize algebraically the defect of stationary lifts for strongly Levi nondegenerate quadrics.

\subsection{Nondefective stationary lifts for strongly Levi nondegenerate quadrics: a special case}

Let $Q=\{\rho=0\}\subset \C^{n+d}$ be a strongly Levi nondegenerate quadric given by \eqref{eqquad}. We consider $b \in \R^d$
such that  $A:=\sum b_jA_j$ is invertible and a stationary lift ${\bm f} \in \mathcal{S}_0(Q)$ of the form \eqref{eqinitlift}. We write explicitly the $(2n+2d)\times (2n+2d)$ matrix 
$\partial \tilde{\rho}({\bm f},\overline{{\bm f}})(\zeta)$. For convenience, we reorder the variables $(z,w,\tilde{z},\tilde{w})$ as  $(w,z,\tilde{z},\tilde{w})$, and so 
$\partial \tilde{\rho}=(\partial_w \tilde{\rho},\partial_z \tilde{\rho},\partial_{\tilde z} \tilde{\rho},\partial_{\tilde{w}} \tilde{\rho})$. We have

\begin{equation}\label{eqG}
\partial \tilde{\rho}\left(\bm{f},\overline{\bm{f}}\right)(\zeta)=\left(\begin{array}{ccccccccccccc}
\frac{1}{2}I_d & \overline{\mathcal{B}(\zeta)} & 0 & 0\\
 0    & \overline{\zeta} A& I_n &   \overline{\mathcal{C}(\zeta)}  \\
  0      & -i\overline{\zeta} A & iI_n  &  i\overline{\mathcal{C}(\zeta)}  \\
0 & 0   &  0  &  i\overline{\zeta} I_d\\
\end{array}\right)
\end{equation}
with 
\begin{equation}\label{eqB}
 \overline{\mathcal{B}(\zeta)}=\left(\begin{matrix}
- \overline{(A_1)^1 h(\zeta)}& \hdots  &- \overline{(A_1)^n h(\zeta)} \\
\vdots &  & \vdots   \\

- \overline{(A_d)^1 h(\zeta)}& \hdots  &- \overline{(A_d)^n h(\zeta)}  \\

\end{matrix}\right)=-(1-\overline{\zeta}) \transp \overline{D}
\end{equation}
and 
\begin{equation}\label{eqC}
\overline{\mathcal{C}(\zeta)}=\left(\begin{matrix}
2 \transp \overline{h} (A_1)_1 & \hdots  & 2 \transp \overline{h} (A_d)_1 \\
\vdots &  & \vdots    \\
2 \transp \overline{h} (A_1)_n & \hdots  & 2 \transp \overline{h} (A_d)_n \\
\end{matrix}\right)=2(1-\overline{\zeta})\overline{D}.
\end{equation}
We recall that $(A_j)_l$ (resp. $(A_j)^l$) denotes, for $j=1,\dots,d$ and $l=1,\ldots,n$, the $l^{\rm th}$ column (resp. the $l^{\rm th}$ row) of $A_j$. We also recall that $D$ is the $n\times d$ matrix whose columns are $A_1V,\ldots,A_dV$.   
We aim to find an algebraic characterization of the nondefectiveness of ${\bm f}$ which is similar to the first point of Lemma \ref{lemnondef}. In this direction, we have the following theorem: 
\begin{theo}\label{theodefli}
Let $Q=\{\rho=0\}\subset \C^{n+d}$ be a strongly Levi nondegenerate  quadric given by \eqref{eqquad}.     
A stationary lift ${\bm f}=(f,\tilde{f})$ of the form \eqref{eqinitlift}  is nondefective if and only if  the $d\times d$ matrix $\Re e [\transp \overline D A^{-1}D]$ is invertible. 
In other words, such a stationary lift is nondefective if and only if $M$ is $\mathfrak{D}$-nondegenerate at $0$. 
Moreover, we have in general 
\begin{equation}\label{eqdeffla}
\bm{{\rm def}}(\bm{f})=\dim {\rm ker} \ \Re e [\transp \overline D A^{-1}D]=d- {\rm rank}\ \Re e [\transp \overline D A^{-1}D].
\end{equation}
\end{theo}

\begin{proof}
Let ${\bm f}$ be a stationary lift of the form \eqref{eqinitlift}.

Assume first  that ${\bm f}$ is defective. We then consider $\tilde{c}: \partial \Delta \to \R^{2n+2d}$ with $\tilde{c}=(\gamma,u,v,\delta)$ with  
$\gamma,\delta: \partial \Delta \to \R^d$ and $u,v: \partial \Delta \to \R^n$ such that the map 
$$\zeta \mapsto \tilde{c}(\zeta)\partial \tilde{\rho}({\bm f},\overline{\bm f})(\zeta)$$ 
defined on $\partial \Delta$ extends holomorphically to $\Delta$. 
Using \eqref{eqG}, \eqref{eqB} and \eqref{eqC}, we compute $\tilde{c}\partial \tilde{\rho}({\bm f},\overline{\bm f})$   and obtain
 \begin{equation}\label{eqdelift}
\Small \left(\frac{1}{2}\gamma(\zeta),-(1-\overline{\zeta})\gamma(\zeta) \transp\overline{D}+\overline{\zeta} (u(\zeta) -i v(\zeta))A ,u(\zeta)+iv(\zeta),
2(u(\zeta)+iv(\zeta))(1-\overline{\zeta})\overline{D}+i\delta(\zeta)\overline{\zeta}\right).
\end{equation}
Note that the holomorphic extension of the first component implies that $\gamma=\gamma_0 \in \R^d$ is constant.  
From the third component, we get directly that $u(\zeta)+iv(\zeta)$ extends as a holomorphic map to $\Delta$.  In particular, we have for $j\geq 1$ 
\begin{equation}\label{eqcoef}
\overline{u_j}+i\overline {v_j}=0
\end{equation}
where $u= \sum_{j\geq 1} \overline{u_j}  \zeta^{-j} + \sum_{j \geq 0} u_j \zeta^j$ and $v=\sum_{j\geq 1} \overline{v_j}\zeta^{-j} + \sum_{j\geq 0} v_j \zeta^j$.  
Looking now at the $\overline{\zeta}$ coefficient of the second  component, we obtain
\begin{equation}\label{eqdeflift}
\gamma_0\  \transp\overline{D}+(u_0-iv_0)A=0,
\end{equation}
while the vanishing of the $\overline{\zeta}^j$ coefficient for $j\geq 2$  gives, using \eqref{eqcoef},
$$(\overline{u_{j-1}}-i\overline{v_{j-1}})A=2\overline{u_{j-1}}A=0.$$
Since $A$ is invertible, we have $u_{j-1}=0$ for $j\geq 2$, and so $u=u_0 \in \R^n$ and $v=v_0 \in \R^n$ are both constant. Back to \eqref{eqdeflift}, we get 
\begin{equation*}
-\gamma_0\ ^{t}\overline{D}=(u_0-iv_0)A
\end{equation*}
and since $A$ is invertible, we obtain
\begin{equation}\label{eqdeflift2}
-\gamma_0\transp\overline{D}A^{-1}=u_0-iv_0.
\end{equation}
 We now look at the fourth component. Note first that since $u$ and $v$ are constant, the map $\delta$ is also constant $\delta=\delta_0 \in \R^d$. Moreover, the $\overline{\zeta}$ coefficient of that fourth component gives
\begin{equation}\label{eqgamma_0}
-2(u_0+iv_0)\overline{D}+i\delta_0=0.
\end{equation}
Using \eqref{eqdeflift2}, we obtain 
\begin{equation*}
\delta_0=2i\gamma_0\transp D\overline{A}^{-1}\overline{D}\in \R^d
\end{equation*}
which implies that  
$$\gamma_0\Re e [\transp \overline{D}A^{-1}D]=0.$$ 
Note finally that  $\gamma_0 \neq 0$, otherwise $u,v$ and $\delta$, and so $\tilde{c}$,  are identically equal to zero.     
So the $d\times d$ matrix $\Re e [\transp \overline{D}A^{-1}D]$ is noninvertible.


\vspace{0.5cm}

Assume now that  $\Re e [\transp \overline D A^{-1}D]$ is noninvertible and let $\gamma_0 \in \R^d\setminus\{0\}$ be such that  
$$\gamma_0\Re e [\transp \overline{D}A^{-1}D]=0.$$ 
So we have $\gamma_0\transp  D \overline{A}^{-1}\overline{D} \in i\R^d.$ 
We then claim that for 
$$\tilde{c}:=\left(\gamma_0,-\Re e \left(\gamma_0\transp\overline{D}A^{-1}\right),\Im m \left(\gamma_0\transp\overline{D}A^{-1}\right),2i\gamma_0\transp D\overline{A}^{-1}\overline{D}\right) 
\in \R^d\times \R^n \times \R^n\times \R^d$$
the map $\zeta \mapsto \tilde{c} \partial \tilde{\rho}({\bm f},\overline{\bm f})(\zeta)$ extends holomorphically to $\Delta$. This follows from  \eqref{eqdelift} and all computations 
performed above. Indeed, the vanishing of all $\overline{\zeta}^j$ coefficients for  $j\geq 2$ is immediate since $\tilde{c}$ is constant, and the vanishing of all 
$\overline{\zeta}$ coefficients follows from \eqref{eqdeflift} and \eqref{eqgamma_0}.  This shows that ${\bm f}$ is defective.  


\vspace{0.5cm}

Finally, it remains to establish \eqref{eqdeffla}. Assume that ${\bm f}$ is defective. As seen above, the map $\tilde{c}$ is  constant in that case. Moreover the linear map 
$$\gamma_0\mapsto \left(\gamma_0,-\Re e \left(\gamma_0\transp\overline{D}A^{-1}\right),\Im m \left(\gamma_0\transp\overline{D}A^{-1}\right),2i\gamma_0\transp D\overline{A}^{-1}\overline{D}\right)$$ defined on  the kernel of $\Re e [\transp \overline D A^{-1}D]$ is an isomorphism onto the space of constant $\tilde{c} \in \R^{2n+2d}$ such that  $\zeta \mapsto \tilde{c}\partial \tilde{\rho}({\bm f},\overline{{\bm f}})(\zeta)$ extends holomorphically to $ \Delta$. This gives directly
 $$\bm{{\rm def}}(\bm{f})=\dim {\rm ker}\Re e [\transp \overline D A^{-1}D].$$
\end{proof}

\begin{remark}\label{remstr}  When $Q\subset \C^{n+d}$ is a strongly pseudoconvex quadric then  
the following statements are equivalent:
\begin{enumerate}[i.]
\item the stationary disc $f=((1-\zeta)V,g)$, $V\in \C^n$,  is nondefective, 
\item the real dimension of ${\rm span}\{A_1V,\ldots,A_dV\}$ is equal to $d$,
\item the $d\times d$ matrix $\Re e [\transp \overline D D]$  is invertible,
\item the stationary lift ${\bm f}=(f,\tilde{f})$ of the form \eqref{eqinitlift} is nondefective.
\end{enumerate}
Indeed, in such case, the linear combination $A$ can be chosen to be the identity $I_n$ and, thus,  $\Re e [\transp \overline D D]$  is invertible if and only if 
${\rm dim}_{\R}{\rm span}\{A_1V,\ldots,A_dV\}= d.$ Moreover, we have 
$${\rm def}(f)=\dim {\rm ker}D_{|\R^d}= \dim {\rm ker} \Re e [\transp \overline D D]=\bm{{\rm def}}(\bm{f}).$$ 
\end{remark}

\begin{example}\label{exdef}
In the context of Example \ref{exline},  Theorem \ref{theodefli} implies that $\bm{{\rm def}}(\bm{f})=1$. Indeed, we have $b=(0,0,1)$, $A=A_3=A^{-1}$, and so
$$\Re e [\transp \overline D A^{-1}D]=\transp \overline D A^{-1}D=
\begin{pmatrix}
0& 0 & 1  \\ 0 &0 & 1 \\ 1 &1 & 0 
\end{pmatrix}$$
which is of rank $2$. In fact, for any invertible real linear combination $A=\sum_{j=1}^3 b_jA_j$, the rank of $\Re e [\transp \overline D A^{-1}D]$ is 2. It follows that all the lifts of the stationary disc $f$ are defective with defect equal to $1$. 
\end{example}

We now make an important connection with the unique jet determination of CR automorphisms of a nondegenerate real submanifold. We proved in \cite{be-me} the following result: 
\begin{theo}[\cite{be-me}] \label{theojet1} 
Let $M\subset \C^{n+d}$ be a $\mathcal{C}^4$   real submanifold given by $\eqref{eqred}$.    
We assume that $M$ is \ $\mathfrak{D}$-nondegenerate at $0$.  Then any germ  at $0$ of  CR automorphism  of $M$ of class $\mathcal{C}^3$   is  uniquely determined by its $2$-jet at $0.$     
\end{theo}

In view of Theorem \ref{theodefli} and Theorem \ref{theojet1}, it is now clear that the existence of a nondefective stationary lift of a special form ensures the unique  $2$-jet  determination of CR automorphisms  of a real submanifold.
\begin{cor}\label{corjet2}  
Let $M\subset \C^{n+d}$ be a $\mathcal{C}^4$ real submanifold given by $\eqref{eqred}$. Consider the model quadric $Q$ associated with $M$, that is the quadric defined by considering only terms of order less than two in $\eqref{eqred}$.       
Assume that there exists a nondefective stationary lift ${\bm f}$ of the form \eqref{eqinitlift} for $Q$. Then any germ  at $0$ of  CR automorphism  of $M$ of class $\mathcal{C}^3$   is  uniquely determined by its $2$-jet at $0.$ 
\end{cor}
We will establish a much more general result in Section 2.5. 
We end this section with a remark related to Proposition \ref{proplift}.

\begin{remark}\label{remconv}

Consider a quadric $Q=\{\rho=0\}\subset \C^{n+d}$ given by \eqref{eqquad} and a defective stationary lift ${\bm f} \in \mathcal{S}_0(Q)$ of the form \eqref{eqinitlift}.
 We know from the proof of Theorem  \ref{theodefli}  that there exists a constant vector
$\tilde{c}_0=(\gamma_0,u_0,v_0,\delta_0) \in \R^{2n+2d}$  such that the map 
$$\zeta \mapsto \tilde{c}_0\partial \tilde{\rho}(\bm{f},\overline{\bm{f}})(\zeta)$$
extends holomorphically to $\Delta$. Note that the term   $\partial_\zeta \tilde{r}$ in the matrix $\partial\tilde{\bm{r}}$ is equal to    $\partial_\zeta \tilde{r}=\transp(0,0,0,-i\overline{\tilde{w}})$. 
The analytic disc $\zeta \mapsto  (\bm{f}(\zeta),\zeta)$ is then defective in case one can find a real-valued map  
$\lambda: \partial \Delta \to \R$ be such that 
$$\zeta \mapsto \left(\tilde{c}_0\partial \tilde{\rho}(\bm{f},\overline{\bm{f}})(\zeta),-\frac{i}{2}\delta_0 \cdot b\overline{\zeta}+\lambda(\zeta)\overline{\zeta}\right)$$
 extends holomorphically to $\Delta$, that is, when the map
 $$\zeta \mapsto \left(-\frac{i}{2}\delta_0 \cdot b+\lambda(\zeta)\right)\overline{\zeta}$$ extends holomorphically to $\Delta$.  It turns out that this can only happen if 
 $\delta_0$ is orthogonal to $b$, in which case one can chose $\lambda$ to be identically equal to zero.

 \end{remark}

\subsection{Nondefective stationary lifts for strongly Levi nondegenerate quadrics: the general case}

We now investigate the defectiveness of an arbitrary stationary lift for a strongly Levi nondegenerate quadric $Q=\{\rho=0\}\subset \C^{n+d}$ given by \eqref{eqquad}. 
In light of Theorem \ref{theodefli} and Corollary \ref{corjet2}, it is natural to expect  the existence of a nondefective stationary lift to be connected with the unique $2$-jet determination of  CR 
automorphisms of $Q$ and its deformations.  In this spirit, we consider the  $1$-jet map at $\zeta=1$ defined on $\mathcal{S}_0(Q)$ 
by  
\begin{equation}\label{eq1jet}
\mathfrak j_{1}:\bm{f} \mapsto (\bm{f}(1),\bm{f}'(1)).
\end{equation}
Since  $\bm{f} \in \mathcal{S}_0(Q)$, we have $\bm{f}(1)=(0,0,0,b/2)$ with $b \in \R^d$ fixed and  such that $\sum_{j=1}^d b_jA_j$ is invertible, and so we identify $\mathfrak j_{1}$ with the derivative map  
$\bm{f} \mapsto \bm{f}'(1)$. Following \cite{be-me2} and due to the explicit expression of stationary lifts \eqref{eqform},  the $1$-jet map 
$$\mathfrak j_{1}: \C^d \times \C^n \to  \C^n \times \R^d \times \R^d$$  can be written as 
\begin{equation*} 
(a,V) \mapsto \left(V, \transp\overline{V}(I-\transp\overline {X})\Psi(A_1)(I-X)V, \ldots,\transp\overline{V}(I-\transp\overline {X})\Psi(A_d)(I-X)V,\Im m a\right),
\end{equation*}
where for $j=1, \dots,d$ 
\begin{equation}\label{eqpsi}
\Psi(A_j):=K_j
\end{equation}
and $K_j$ is defined in \eqref{eqK_j}. 
Here, $\Psi$ is a self-map of the the space of $n \times n$ matrices with complex entries. 
We also recall that $X$, which depends on both $a$ and $b$, is the unique $n\times n$ matrix solution of \eqref{eqX} with  $\|X\|<1$.
We are now in a position to state our main theorem. 

\begin{theo}\label{theonondef}
Let  $Q=\{\rho=0\}\subset \C^{n+d}$ be a strongly Levi nondegenerate quadric given by \eqref{eqquad}. Consider a stationary lift ${\bm f} \in \mathcal{S}_0(Q)$ of the form \eqref{eqform}, that is, parametrized by $(a,V) \in \C^d\times \C^n$. Then  ${\bm f}$   is nondefective if and only if  the $1$-jet map $\mathfrak j_{1}$ is a local diffeomorphism at the stationary lift in $\mathcal{S}_0(Q)$ parametrized by 
$(a,(I-X)V)$.  
 \end{theo}

\begin{remark} 
In the context of Theorem \ref{theonondef}, if ${\bm f}$ is nondefective then, due to its analyticity, the $1$-jet map $\mathfrak j_{1}$ is a local diffeomorphism outside of a thin set. 
\end{remark}  
 
Before moving to the proof of Theorem \ref{theonondef}, we need some preparatory setup. We recall a key lemma: 
\begin{lemma}\label{lemjet}\cite{be-me2,be-me-ko}  Let $Q$ be a strongly Levi nondegenerate quadric and let $b \in \R^d$. We denote by  $X_{\Re e  a_s}$ the partial derivative of $X$ with respect to $\Re e  a_s$, $s=1,\ldots,d$. 
We  have:
 \begin{enumerate}[i.]
\item  The $1$-jet map $\mathfrak j_{1}$ is a local diffeomorphism at $(a,V) \in \C^d\times \C^n$ if and only if 
the  $d\times d$ matrix 
\begin{equation}\label{eqjet}
\Re e   \left(\sum_{r=0}^{\infty} \transp\overline{V}\left(I-\transp \overline {X}\right)^2 {\transp \overline {X}}^r\Psi(A_j)X_{\Re e  a_s}X^rV\right)_{j,s}
\end{equation}
 is invertible.
\item For $s=1,\ldots,d$, we have
 \begin{equation}\label{eqbre}
 X_{\Re e a_s}= \sum_{r=0}^\infty (-1)^{r+1}\left((I+A^{-1}P X)^{-1}A^{-1}P\right)^r(I+A^{-1}P X)^{-1}A^{-1} A_s(I-X)^2X^r.
  \end{equation}
\end{enumerate}
\end{lemma}

In relation with the quadratic equation \eqref{eqX}, according to \cite{be-me2} (see also \cite{tu}), the invertible Hermitian $n\times n$  matrix $B:=A+ PX$ is such that 
 for all $\zeta \in \partial \Delta$
\begin{equation}\label{eqfact}
\overline{\zeta}P+A+\zeta\transp \overline{P}=(I-\overline {\zeta}\  \transp{\overline X})B(I- {\zeta}X).
\end{equation}    
This factorization generalizes a   result due to Lempert (Theorem B p. 442 in \cite{le}) for positive definite matrices. 
We can then rewrite Equation \eqref{eqbre} in a more convenient form.
\begin{lemma}\label{lemdifX}
We have 
\begin{equation*}
 X_{\Re e a_s}   =   B^{-1}\Psi(A_s)(I-X)^2
\end{equation*}
where the map $\Psi$ is defined in \eqref{eqpsi}. 
\end{lemma}
 \begin{proof}
 Since $B=A+ PX=A(I+A^{-1}PX)$ and $BX=-\transp\overline {P}$, we write for $s=1,\ldots,d$
 \begin{eqnarray*}
 X_{\Re e a_s}  & = &  \sum_{r=0}^\infty (-1)^{r+1}\left(B^{-1} P\right)^rB^{-1} A_s(I-X)^2X^r\\
 & =& - \sum_{r=0}^\infty B^{-1} \transp \overline{X}^{r-1}PB^{-1} A_sX^r (I-X)^2\\
 &=&  \sum_{r=0}^\infty B^{-1} \transp \overline{X}^{r} A_sX^r (I-X)^2=B^{-1}\Psi(A_s)(I-X)^2. 
  \end{eqnarray*}
\end{proof}
Now following \eqref{eqjet}, we define the  $d\times d$ matrix 
$$\mathcal{M}_1:= \left(\sum_{r=0}^{\infty} \transp\overline{V}\left(I-\transp \overline {X}\right)^2 {\transp \overline {X}}^r\Psi(A_j)X_{\Re e a_s}X^rV\right)_{j,s}$$
and, using Lemma \eqref{lemdifX}, we observe that
 \begin{eqnarray}\label{eqm1}
\notag \mathcal{M}_1 &=&\left(\transp\overline{V}\left(I-\transp \overline {X}\right)^2  \left(\sum_{r=0}^{\infty}  {\transp \overline {X}}^r\Psi(A_j)B^{-1}\Psi(A_s)X^r\right)(I-X)^2V\right)_{j,s} \\
&=& \left(\transp\overline{V}\left(I-\transp \overline {X}\right)^2  \Psi\left(\Psi(A_j)B^{-1}\Psi(A_s)\right)(I-X)^2V\right)_{j,s}.
  \end{eqnarray}

We can now prove our main theorem.
\begin{proof}[Proof of Theorem \ref{theonondef}]
Consider a  stationary lift  ${\bm f}=(h,g,\tilde{h},\tilde{g}) \in \mathcal{S}_0(Q)$  of the form \eqref{eqform}, that is, parametrized by $(a,V) \in \C^d\setminus\{0\} \times \C^n$. The case $a=0$ was already covered in Theorem \ref{theodefli}.  
We first write the $(2n+2d)\times (2n+2d)$  matrix $\partial \tilde{r}({\bm f},\overline{{\bm f}})(\zeta)$
\begin{equation*}
\partial \tilde{r}(\bm{f},\overline{\bm{f}})(\zeta)=\left(\begin{array}{ccccccccccccc}
\frac{1}{2}I_d & \overline{\mathcal{B}(\zeta)} & 0 & 0\\
 0    & \transp \overline{P}+\overline{\zeta} A+\overline{\zeta}^2P& I_n &   \overline{\mathcal{C}(\zeta)}  \\
  0      & -i(\transp \overline{P}+\overline{\zeta} A+\overline{\zeta}^2P) & iI_n  &  i\overline{\mathcal{C}(\zeta)}  \\
0 & 0   &  0  &  i\overline{\zeta} I_d\\
\end{array}\right),
\end{equation*}
with 
\begin{equation*}
 \overline{\mathcal{B}(\zeta)}=\left(\begin{matrix}
- \overline{(A_1)^1 h(\zeta)}& \hdots  &- \overline{(A_1)^n h(\zeta)} \\
\vdots &  & \vdots  \\

- \overline{(A_d)^1 h(\zeta)}& \hdots  &- \overline{(A_d)^n h(\zeta)}  \\

\end{matrix}\right)=-\left(1-\overline{\zeta}\right) \transp \overline{\mathcal{D}}(\zeta)
\end{equation*}
and 
$$\overline{\mathcal{C}(\zeta)}=\left(\begin{matrix}
2 \transp \overline{h(\zeta)} (A_1)_1 & \hdots  & 2 \transp \overline{h(\zeta)} (A_d)_1 \\
\vdots &  & \vdots   \\
2 \transp \overline{h(\zeta)} (A_1)_n & \hdots  & 2 \transp \overline{h(\zeta)} (A_d)_n \\
\end{matrix}\right)=2\left(1-\overline{\zeta}\right)\overline{\mathcal{D}}(\zeta).$$
Here, $\mathcal{D}(\zeta)$ is the following $n\times d$ matrix 
$$\mathcal{D}(\zeta)=\sum_{r=0}^{\infty}\zeta^rD_{r}$$
 where $D_r$ is the matrix whose columns are $A_1X^r V,\ldots,A_dX^r V$. 
%
For $\tilde{c}=(\gamma,u,v,\delta): \partial \Delta \to \R^{d} \times \R^n \times \R^n \times \R^d$, the four components of the map $\zeta \mapsto \tilde{c}(\zeta)\partial \tilde{r}({\bm f},\overline{{\bm f}})(\zeta)$ are then respectively
\begin{equation}\label{eqdelift2}
\begin{cases}
\displaystyle \frac{1}{2}\gamma(\zeta)\\
\\
-\left(1-\overline{\zeta}\right)\gamma(\zeta) \transp\overline{\mathcal{D}}(\zeta)+(u(\zeta) -i v(\zeta))\left(\transp \overline{P}+\overline{\zeta} A+\overline{\zeta}^2P\right)\\
\\
u(\zeta)+iv(\zeta)\\
\\
2\left(1-\overline{\zeta}\right)(u(\zeta)+iv(\zeta))\overline{\mathcal{D}}(\zeta)+i\delta(\zeta)\overline{\zeta}
 \end{cases}
\end{equation}

\vspace{0.5cm}

Assume now that the stationary lift  ${\bm f}=(h,g,\tilde{h},\tilde{g})$ is defective, that is, for $\tilde{c}$ not identically equal to zero, the four components \eqref{eqdelift2} extend holomorphically to $\Delta$. We first note that $\gamma(\zeta)\equiv \gamma_0$ is constant 
and that $u(\zeta)+iv(\zeta)$ is holomorphic on $\Delta$. Writing $u= \sum_{j\geq 1} \overline{u_j}  \zeta^{-j} + \sum_{j \geq 0} u_j \zeta^j$ and $v=\sum_{j\geq 1} 
\overline{v_j}\zeta^{-j} + \sum_{j\geq 0} v_j \zeta^j$, we have once more  for $j\geq 1$,  $\overline{u_j}+i\overline {v_j}=0.$
Moreover, using \eqref{eqfact}, we write the second component as 
\begin{equation}\label{eqsec}
-\left(1-\overline{\zeta}\right)\gamma_0\transp\overline{\mathcal{D}}(\zeta)+(u(\zeta) -i v(\zeta))(I-\overline{\zeta}\transp\overline{X})B\left(\overline{\zeta}- X\right)
\end{equation}
and we observe that it is constant as well and equal to 
\begin{equation*}
-\gamma_0\transp\overline{D_0}-(u_0-i v_0)BX.
\end{equation*}
We right-multiply the map \eqref{eqsec}   by $\left(\overline{\zeta}-X\right)^{-1}=\zeta(I-\zeta X)^{-1}$ and obtain  
\begin{equation}\label{eqab}
-\left(1-\overline{\zeta}\right)\gamma_0\transp\overline{\mathcal{D}}(\zeta)\zeta(I-\zeta X)^{-1}+(u(\zeta) -i v(\zeta))\left(I-\overline{\zeta}\transp\overline{X}\right)B
\end{equation}
which also extends holomorphically to $\Delta$.  We then study for $r\geq 0$ the $\overline{\zeta}^r$  coefficient of 
$$\left(1-\overline{\zeta}\right)\gamma_0\transp\overline{\mathcal{D}}(\zeta)\zeta(I-\zeta X)^{-1}=\gamma_0\left(\transp \overline{D_0}+\sum_{r\geq 1} \overline{\zeta}^r\transp \overline{(D_r-D_{r-1})}\right)\sum_{r\geq 1} \zeta^{r}X^{r-1}.$$
Its constant coefficient is equal to
$$\gamma_0\sum_{k\geq 1} \transp \overline{(D_k-D_{k-1})}X^{k-1}$$
and, in general, its  $\overline{\zeta}^r$ coefficient for $r\geq 1$ is
$$\gamma_0\sum_{k\geq r+1} \transp \overline{(D_k-D_{k-1})}X^{k-r-1}.$$
It follows that the constant, $\overline{\zeta}$ and $\overline{\zeta}^r$ for  $r\geq 2$,  coefficients of the map \eqref{eqab} are respectively
$$-\gamma_0\sum_{k\geq 1} \transp \overline{(D_k-D_{k-1})}X^{k-1}+(u_0-iv_0)B =0$$
$$-\gamma_0\sum_{k\geq 2} \transp \overline{(D_k-D_{k-1})}X^{k-2}+(-(u_0-iv_0)\transp\overline{X}+2\overline{u_1})B=0$$
and
$$-\gamma_0\sum_{k\geq r+1} \transp \overline{(D_k-D_{k-1})}X^{k-r-1}+(-2\overline{a_{r-1}}\transp \overline{X}+2\overline{u_r})B=0.$$
From this triangular system of linear equations,  we obtain the following 
%
%
%
%

\begin{eqnarray*}
u_0-iv_0&=&\gamma_0\sum_{k\geq 1} \transp \overline{(D_k-D_{k-1})}X^{k-1}B^{-1}\\
\\
2\overline{u_1}&=&(u_0-iv_0)\transp\overline{X} + \gamma_0\sum_{k\geq 2} \transp \overline{(D_k-D_{k-1})}X^{k-2}B^{-1}\\
&=& \gamma_0\sum_{k\geq 1} \transp \overline{(D_k-D_{k-1})}X^{k-1}B^{-1}\transp\overline{X} +\gamma_0\sum_{k\geq 2} \transp \overline{(D_k-D_{k-1})}X^{k-2}B^{-1}\\
& =&  \gamma_0 \sum_{\ell=0}^1 \sum_{k\geq \ell+1} \transp \overline{(D_k-D_{k-1})}X^{k-\ell-1}B^{-1}\transp\overline{X}^{1-\ell}\\
\\
 2\overline{u_r}&=& \gamma_0\sum_{k\geq r+1} \transp \overline{(D_k-D_{k-1})}X^{k-r-1}+2\overline{a_{r-1}}\transp \overline{X}\\
& =& \gamma_0 \sum_{\ell=0}^r \sum_{k\geq \ell+1} \transp \overline{(D_k-D_{k-1})}X^{k-\ell-1}B^{-1}\transp\overline{X}^{r-\ell}.\\
\end{eqnarray*}
%
Thus, we have for $r\geq 1$
\begin{equation}\label{eqa}
\begin{cases}
u_0+iv_0=\gamma_0\sum_{k\geq 1} \transp(D_k-D_{k-1})\overline{X}^{k-1}\overline{B}^{-1}\\
\\
 2u_r= \gamma_0 \sum_{\ell=0}^r \sum_{k\geq \ell+1} \transp(D_k-D_{k-1})\overline{X}^{k-\ell-1}\overline{B}^{-1}\transp X^{r-\ell}.\\
\end{cases}
\end{equation}
We now look at the fourth component of \eqref{eqdelift2} whose $\overline{\zeta}$ coefficient gives
$$2(u_0+iv_0)\left(\overline{D_1}-\overline{D_0}\right)+4\sum_{r\geq 1}u_r\left(\overline{D_{r+1}}-\overline{D_r}\right)+i\delta_0=0$$
which implies that 
$$(u_0+iv_0)\left(\overline{D_1}-\overline{D_0}\right)+2\sum_{r\geq 1}u_r\left(\overline{D_{r+1}}-\overline{D_r}\right)\in i\R^d.$$
Using \eqref{eqa}, this last expression can be written as
%
\begin{equation}\label{eqfourth}
\gamma_0\sum_{r\geq 0}\left[\sum_{\ell=0}^r \sum_{k\geq \ell+1} \transp(D_k-D_{k-1})\overline{X}^{k-\ell-1}\overline{B}^{-1}\transp X^{r-\ell}\right]\left(\overline{D_{r+1}}-\overline{D_r}\right) \in  i\R^d\\
\end{equation}
Recall that $D_r$ is the  $n\times d$ matrix whose $j^{th}$ column is 
$A_jX^r V$. So we have  for $r\geq 1$
$$A_jX^r V-A_jX^{r-1} V=A_jX^{r-1}\underbrace{(X-I)V}_{V'}.$$   
We then write \eqref{eqfourth} as 
$$\gamma_0\left(\sum_{r\geq 0}\sum_{\ell=0}^r \sum_{k\geq \ell+1}  \transp V' \transp X^{k-1}\overline{A_j}\overline{X}^{k-\ell-1}
\overline{B}^{-1}\transp X^{r-\ell}\overline{A_sX^rV'}\right)_{j,s} \in  i\R^d.$$
The real part of the above expression is then identically equal to zero. This implies that the {\it real part} of the following $d\times d$ matrix
$$\mathcal{M}_2:=\left(\sum_{r\geq 0}\sum_{\ell=0}^r \sum_{k\geq \ell+1}  \transp \overline{V'} \transp \overline{X}^{k-1}A_j X^{k-\ell-1}B^{-1}\transp \overline{X}^{r-\ell}A_sX^rV'\right)_{j,s}$$
is noninvertible. To conclude, we write $\mathcal{M}_2$ in a form similar to  \eqref{eqm1}. It is important to note that this is possible since all below matrix series converge uniformly in view of  $\|X\|<1$.
 \begin{eqnarray*}
\mathcal{M}_2& = &  \transp \overline{V'}  \sum_{r= 0}^\infty\sum_{\ell=0}^r \sum_{k= 0}^{\infty} \transp \overline{X}^{k+\ell}A_j X^{k}B^{-1}\transp \overline{X}^{r-\ell}A_sX^rV'\\
& = & \transp \overline{V'}  \sum_{r= 0}^\infty \sum_{\ell=0}^r \transp \overline{X}^{\ell} \sum_{k= 0}^{\infty} \transp \overline{X}^{k}A_j X^{k}B^{-1}\transp \overline{X}^{r-\ell}A_sX^rV'\\
& = & \transp \overline{V'}  \sum_{r= 0}^\infty \sum_{\ell=0}^r \transp \overline{X}^{\ell}\Psi(A_j)B^{-1}\transp \overline{X}^{r-\ell}A_sX^{r-\ell}X^\ell V'\\
& = & \transp \overline{V'}  \sum_{r= 0}^\infty \sum_{\ell=0}^r \transp \overline{X}^{r-\ell}\Psi(A_j)B^{-1}\transp \overline{X}^{\ell}A_sX^{\ell}X^{r-\ell} V'\\
& = & \transp \overline{V'}  \sum_{\ell=0}^\infty  \sum_{r= \ell}^\infty \transp \overline{X}^{r-\ell}\Psi(A_j)B^{-1}\transp \overline{X}^{\ell}A_sX^{\ell}X^{r-\ell} V'\\
& = & \transp \overline{V'}  \sum_{\ell=0}^\infty \Psi\left(\Psi(A_j)B^{-1}\transp \overline{X}^{\ell}A_sX^{\ell}\right) V'\\
& =& \transp\overline{V'}\Psi\left[\Psi(A_j)B^{-1}\Psi(A_s)\right]V'.
\end{eqnarray*}
Since $\Re e \mathcal{M}_2$ is noninvertible, it follows from Equation \eqref{eqm1} that the $1$-jet map $\mathfrak j_{1}$ is not a local diffeomorphism at $(a,(I-X)V).$

 \vspace{0.5cm}
 
 We now prove the converse. Assume that the $1$-jet map $\mathfrak j_{1}$ is not a local diffeomorphism at the stationary lift parametrized by  $(a,(I-X)V) \in \C^d\setminus\{0\} \times \C^n.$ Based on the above, there exists 
 $\gamma_0\in \R^d\setminus\{0\} \cap {\rm ker} (\Re e \mathcal{M}_2)$, that is, such that  $\gamma_0 \mathcal{M}_2 \in i\R^d$. We then consider
 $$\tilde{c}(\zeta)=(\gamma_0,u(\zeta),v(\zeta),\delta(\zeta))$$
 where the map $u+iv$ is holomorphic on $\Delta$ with coefficients defined by Equation \eqref{eqa} and where the coefficients of the map $\delta$ are such that $\delta_0=-2\gamma_0 \mathcal{M}_2 \in i\R^d$ and for $r\geq 1$  
 $$\delta_r=2i(u_0+iv_0)\left(\overline{D_r}-\overline{D_{r-1}}\right)+4i\sum_{k\geq 1}u_k\left(\overline{D_{k+r}}-\overline{D_{k+r-1}}\right) \in \C^d.$$
It follows that the map $\zeta \mapsto \tilde{c}(\zeta)\partial \tilde{r}({\bm f},\overline{{\bm f}})(\zeta)$ extends holomorphically to $\Delta$. The  stationary lift  ${\bm f} \in \mathcal{S}_0(Q)$ parametrized by $(a,V)$  is then defective.

\end{proof}

As a direct corollary, we obtain 
\begin{cor}\label{cordef}
Let ${\bm f}=(f,\tilde{f})$ be a stationary lift of the form \eqref{eqform}. Its defect is given by 
$$\bm{{\rm def}}(\bm{f})=\dim {\rm ker}\Re e \mathcal{M}_2.\ \ $$
\end{cor}

We now make an important observation in the strongly pseudoconvex case and which echoes Remark \ref{remstr}. 
\begin{cor}\label{corstrong}
Assume that $Q \subset \C^{n+d}$ is a strongly pseudoconvex quadric. Let $f$ be a stationary disc of the form \eqref{eqform}. Then 
 the  disc $f$ is nondefective if and only if any of its lift ${\bm f}=(f,\tilde{f})$ is nondefective. 
\end{cor}

\begin{proof}
This result follows from two important points. First, if $f=(h,g)$ is nondefective then the quadric $Q$ is stationary minimal at $0$ for $(a,b-a-\overline{a},h'(0))=(a,b-a-\overline{a },(X-I)V)$  (see Proposition 3.3 \cite{be-me2} or also \cite{tu}). Moreover, we know from Theorem 4.3 \cite{be-me2} that if the quadric $Q$ is stationary minimal at $0$ for 
$(a,b-a-\overline{a},(X-I)V)$ then the $1$-jet map $\mathfrak{j}_1$ is a local diffeomorphism at $(a,(X-I)V)$. The claim then follows from Theorem \ref{theonondef}. 
\end{proof}
It is important to point out that, 
when $Q \subset \C^{n+d}$ is a strongly pseudoconvex quadric,  Tumanov proved in \cite{tu3} the 
existence of a nondefective stationary disc attached to $Q$. We illustrate Tumanov's proof in the particular case where $Q$ is of maximal codimension $d=n^2$.

\begin{example}\label{exd=n^2}
Let $Q$ be a quadric of maximal codimension $d=n^2$ satisfying $\aaa$. As pointed out in \cite{be-me4}, $Q$ is  strongly pseudoconvex. 
After a  holomorphic linear change of coordinates, $Q$ can be written as	
\begin{equation*}
\Re e  w_{j,k}=\transp \overline{z} A_{j,k} z \ \ 1\leq j,k \leq n
\end{equation*}
where $A_{j,k}, 1\leq j,k \leq n$ is the canonical basis of the space of Hermitian $n\times n$ matrices. In particular, we have $\Re e  w_{j,j}=|z_j|^2$.
Now for  $V=\transp (1, 1,  \ldots,1) \in \C^n$ and $a=(a_{j,k}) \in \R^{n^2}$ with  $a_{j,k}= 0$ whenever $j\ne k$ and with $a_{j,j}$ pairwise distinct,  the matrix $P$ from Equation \eqref{eqX} is diagonal and so the matrices $A_{j,k}$'s restricted to the  space 
$${\rm span}_\R \{V, PV, {P^2}V, \ldots, {P^k}V, \ldots\}=
{\rm span}_\R \left\{\begin{pmatrix}
	1\\
	\vdots \\
	1\\ 
\end{pmatrix}, \begin{pmatrix}
	a_{1,1}\\
	\vdots \\
	a_{n,n}\\ 
\end{pmatrix}, \begin{pmatrix}
	a_{1,1}^2\\
	\vdots \\
	a_{n,n}^2\\ 
\end{pmatrix}, \ldots, \begin{pmatrix}
	a_{1,1}^k\\
	\vdots \\
	a_{n,n}^k\\ 
\end{pmatrix}, \ldots\right\}$$ 
are $\R$-linearly independent by the invertibility of  the Vandermonde matrix. This implies $Q$ is stationary minimal for $(a,b-2a,V)$ where $b=(b_{j,k}) \in \R^{n^2}$ with 
$b_{j,k}= 0$ whenever $j\ne k$ and $b_{j,j}=1.$ In particular by Lemma \ref{lemstatmin}, the stationary disc parametrized by $(a,V)$ is nondefective. 

\end{example}

To conclude this section, we state an important consequence of  Theorem \ref{theonondef} which follows directly from Proposition 3.6 \cite{be-bl-me} and which
allows us to treat small perturbations of the quadric $Q$. 
\begin{cor}\label{cornondef}
Let $Q=\{\rho=0\}\subset \C^{n+d}$ be a strongly Levi nondegenerate quadric  the form \eqref{eqquad} and let $\bm{f_0}  \in \mathcal{S}_0(Q)$ be a 
nondefective stationary lift of the form \eqref{eqform} parametrized by $(a,V)$ and let $\bm{f_1}  \in \mathcal{S}_0(Q)$ be the stationary lift parametrized by $(a,(I-X)V)$. 
Then there exist a 
neighborhood $U$ of $\rho$ in $(\mathcal{C}^4(\B))^{d}$ and a positive number $\varepsilon>0$ such that for all $r\in U$, the $1$-jet map $\mathfrak j_{1}$ 
$$\mathfrak j_{1}: \{\bm{f} \in \mathcal{S}_0(\{r=0\})\ | \ \|\bm{f}-\bm{f_1}\|_{\mathcal A^{1,\alpha}_{0}}<\varepsilon\} \to \C^{n+d}$$
is a local diffeomorphism at some lift ${\bm f}$, provided ${\bm f}$ is close enough to $\bm{f_1}$.  
\end{cor}

\subsection{Jet determination of CR automorphisms }

Let $k$ be a positive integer.   
Let $M\subset \C^{n+d}$ be a $\mathcal{C}^4$ generic  real submanifold and let $p \in M$. We denote by $Aut^k(M,p)$ the set of germs at $p$  of CR automorphisms $F$ of $M$ of class $\mathcal{C}^k$; in particular we have $F(p)=p$ and $F(M)\subset M$. Combining Theorem \ref{theonondef} and several results from \cite{be-me}, we obtain the following theorem. 
\begin{theo}\label{theojet}
Let $M\subset \C^{n+d}$ be a $\mathcal{C}^4$ generic strongly Levi nondegenerate real submanifold given by \eqref{eqred}. Assume that its model quadric $Q$ admits a nondefective stationary lift ${\bm f}$ of the form \eqref{eqform}. Then elements of $Aut^3(M,0)$ are uniquely determined by their $2$-jet at $0$.     
\end{theo}
As already pointed out in the introduction, this theorem unifies the two previously known results on  $2$-jet determination \cite{tu3, be-me, me, gr-me,be-me-ko}. Its hypothesis is purely geometric and is free of codimension restrictions.

In Proposition 4.4 in \cite{be-me} (see also \cite{tu, sc-tu} in the strongly pseudoconvex case), we observed that the $1$-jet determination at $\zeta=1$ of stationary lifts is a key property that leads to the $2$-jet determination of germs at $0$ of CR automorphisms of $M$. Indeed, the fact that the $1$-jet map $\mathfrak j_{1}$ defined in \eqref{eq1jet} is a local diffeomorphism at some stationary lift ${\bm f} \in \mathcal{S}_0(Q)$ allows us to cover an open set in the conormal bundle. 
\begin{prop}\cite{be-me}\label{propfill3}
Let $Q=\{\rho=0\}\subset \C^{n+d}$ be a strongly Levi nondegenerate quadric  given by \eqref{eqquad} and let $\bm{f_0}  \in \mathcal{S}_0(Q)$ be a 
stationary lift  of the form \eqref{eqform}. Assume that the  map $\mathfrak j_{1}$ is a local diffeomorphism at ${\bm f_0}$. Then there exist a
neighborhood $U$ of $\rho$ in $(\mathcal{C}^4(\B))^{d}$ and a positive number $\varepsilon>0$ such that for all $r\in U$, the 
image of the map
$$\psi:\{\bm{f} \in \mathcal{S}_0(\{r=0\})\ | \ \|\bm{f}-\bm{f_0}\|_{\mathcal A^{1,\alpha}_{0}}<\varepsilon\} \times [0,\varepsilon) \to N^*(\{r=0\})$$
defined by $$\psi({\bm f},\theta)=e^{-i\theta}.{\bm f}(e^{i\theta})=(f(e^{i\theta}),e^{-i\theta}\tilde{f}(e^{i\theta}))$$
covers an open set in the conormal bundle $N^*(\{r=0\})$.

\end{prop}
We point out that, in \cite{be-me}, this proposition is stated for a lift $\bm{f_0}  \in \mathcal{S}_0(Q)$ parametrized $(0,V)$. However the fact $a=0$ is never used in the proof and the main property that the disc must satisfy is that the differential at $\bm{f_0}$ of the  $1$-jet map $\mathfrak j_{1}$ is an isomorphism. Moreover, we already pointed out in Remark \ref{remman}, that 
$$\{\bm{f} \in \mathcal{S}_0(\{r=0\})\ | \ \|\bm{f}-\bm{f_0}\|_{\mathcal A^{1,\alpha}_{0}}<\varepsilon\}$$ is a $\mathcal{C}^1$ real submanifold of dimension $2n+2d$.  

\vspace{0.5cm}

The proof of Theorem \ref{theojet} is standard and  identical to the one of Theorem in 4.7 \cite{be-me} (see also the proof of  Theorem 4.1 in \cite{be-bl-me}). 
For sake of clarity, we present its main points. Assume that  $F$ is a germ at $0$ of  CR automorphism of class $\mathcal{C}^3$ of $M$ fixing $p=0$ and with 
a trivial $2$-jet at $0$. Consider the model quadric $Q$ associated with $M$ and let  $\bm{f_0} \in \mathcal{S}_0(Q)$  be a nondefective stationary lift of the form 
(\ref{eqform}) parametrized by $(a,V)$, and let $\bm{f_1} \in \mathcal{S}_0(Q)$  be the stationary lift parametrized by $(a,(I-X)V)$. We consider the open set 
$\mathcal{O} \subset N^*(M)$ obtained in Proposition \eqref{propfill3}. Let $q \in \mathcal{O}$  and consider a stationary lift $\bm{f} \in \mathcal{S}_0(M)$ for 
$M$ close to $\bm{f_1}$ and such that   $\psi({\bm f},\theta)=q$ for some small $\theta$. Since $F$ has a trivial $2$-jet at $0$, the discs $\bm{f}$  and its push forward
$$F_{*}\bm{f}(\zeta)=\left(F\circ f(\zeta), \tilde{f}(\zeta)\left(d_{f(\zeta)}F\right)^{-1}\right)$$ 
have the same $1$-jet at $\zeta=1$ and 
therefore coincide by Corollary \ref{cornondef}.  It follows that the lift to the cotangent bundle of $F$,  still 
denoted by $F$ satisfies $F(q)=q$. Thus, the map $F$ is the identity on the open set $\mathcal{O} \subset N^*(M)$ and so is identically equal to the identity.

\vskip 1cm

{\small
\noindent Florian Bertrand\\
Department of Mathematics\\\
Center for Advanced Mathematical Sciences\\
American University of Beirut, Beirut, Lebanon\\{\sl E-mail address}: fb31@aub.edu.lb\\

\noindent Francine Meylan \\
Department of Mathematics\\
University of Fribourg, CH 1700 Perolles, Fribourg\\
{\sl E-mail address}: francine.meylan@unifr.ch\\
}

\end{document}